\documentclass[12pt, reqno]{amsart}
\usepackage{amscd,bbold,amsbsy}
\usepackage{amssymb}
\usepackage{amsmath}
\usepackage[all]{xypic}
\usepackage{graphics}
\def\exterior{{\mathchoice{\mbox{\large\boldmath$\wedge$}}
{\mbox{\large\boldmath$ \wedge$}}
{\mbox{\small\boldmath$\wedge$}}
{\mbox{\scriptsize\boldmath$\wedge$}} }}

\def\sym{{\mathchoice{\mbox{\sf S\hspace{0.1ex}}}
{\mbox{\sf S\hspace{0.1ex}}}
{\mbox{\scriptsize\sf S\hspace{0.1ex}}}
{\mbox{\tiny\sf S\hspace{0.1ex}}} }}

\newsymbol\injarrow 131A
\newsymbol\surarrow 1310

\def\hi{\hat{i}}

\def\hj{\hat{j}}

\def\hk{\hat{k}}
\def\hl{\hat{l}}
\def\hp{\hat{p}}
\def\hq{\hat{q}}

\def\hs{\hat{s}}

\def\db{\mathrm{db}} 

\def\ev{\mathrm{ev}}
\def\id{\mathrm{id}}

\newtheorem{thm}{Theorem}[section]

\newtheorem{lem}[thm]{Lemma}

\newtheorem{pro}[thm]{Proposition}
\newtheorem{rmk}[thm]{Remark}
 
\newcommand{\bb}{\begin{equation}}
\newcommand{\eqbb}{\begin{equation}}
\def\ee{\end{equation}}
\def\eqee{\end{equation}}
 
\def\Ker{\mbox{\rm Ker}}

\def\ch{\mbox{\rm ch}}

\def\Im{\mbox{\rm Im}}

\def\End{\mbox{\rm End\hspace{0.1ex}}}

\def\H{{\mathcal H}}

\def\Id{\id}

\def\secteqn{
\let\sectio\section%
\renewcommand{\section}{\sectioneqn\sectio }%
\newcommand{\sectioneqn}{\setcounter{equation}{0}%
 \renewcommand{\theequation}{\arabic{section}.\arabic{equation}}}} 
\title[Irreducible Representations of $GL_q(3|1)$]{Construction of Irreducible Representations of the Quantum  Super Group  $GL_q(3|1)$}

\author[Dung]{NGUY\^EN Thi Phuong Dung} 
  \address[Dung]{Border Academy, Son Tay, Hanoi, Vietnam }
 \email{phuongdung72@yahoo.com}
 
 \author[Hai]{Ph\`ung H\^o Hai}
 \address[Hai]{Institute of Mathematics, Hanoi, Vietnam}
 \email{phung@math.ac.vn}
 
  \author[Hung]{NGUY\^EN Huy Hung}
 \address[Hung]{ Hanoi Pedagogical University II, Vinhphuc, Vietnam}
 \email{hungp1999@yahoo.com}
\thanks{The work is supported to the authors by NAFOSTED under grant no. 101.01.16.9.}
\begin{document}

 \maketitle \begin{abstract} In this note, we  construct  all irreducible representations of the quantum general linear super group $GL_q(3|1)$ using the double Koszul complex.
\end{abstract}
 \bibliographystyle{plain}

\section{Introduction}
A quantum general linear super group is understood as a Hopf super algebra determined in terms of a Hecke symmetry $R$ on a super vector space $V$ of finite dimension. A representation of such a quantum group is nothing but a comodule on the corresponding Hopf super algebra.  

The main invariant of a Hecke symmetry is its birank. 
It is shown in \cite{phhai05} that the category of representations of this quantum group is uniquely determined up to braided monoidal equivalence by the birank of the Hecke symmetry $R$, provided that the quantum parameter $q$ is not a root of unity of order larger than 1.  Therefore, the quantum general linear super group associated to a Hecke symmetry of birank $(r,s)$ is denoted simply by $GL_q(r|s)$.

An explicit construction of irreducible representations, i.e. simple comodules over the associated Hopf super algebra, is however not known. Actually, such a construction is not known even in the classical situation of the Lie super algebras $\mathfrak{gl}(m|n)$. The difficulty lies in the so called atypical representations.

Some particular cases of lower biranks $(1|1)$ and $(2|1)$ are treated in \cite{phhai01,dh2}. Recently, an explicit construction of irreducible representations of $\mathfrak{gl}(3|1)$ was obtained in \cite{Dung1} using the so called double Koszul complex. In this work, this construction will be extended to the case of quantum general linear super group $GL_q(3|1)$. To show that the representations obtained are indeed irreducible and furnish all irreducible representations we use a result of \cite{zhang06} on the perfect paring between $GL_q(r|s)$ and $\mathcal U_q(\mathfrak{gl}(r|s))$ as well as the character formula for these representations.

\section{The quantum general linear supergroup}\label{sect1}

Let $V$ be a super vector space of finite dimension over $\mathbb k$, an algebraically closed field of characteristic zero. Fix a homogeneous basis $x_1,x_2,\ldots, x_d$ of $V$. We shall denote the parity of the basis element $x_i$ by $\hi$. An even operator $R$ on $V\otimes V$ can be given by a matrix $R_{ij}^{kl}$: 
$$R(x_i\otimes x_j)=x_k\otimes x_lR^{kl}_{ij}$$ $R$ is an even operator implies that the matrix elements $R^{ij}_{kl}$ are zero, except for those with $\hi+\hj=\hk+\hl$.
 $R$ is called {\it Hecke symmetry} if the following conditions are satisfied:
\begin{itemize}\item[i)] $R$ satisfies the Yang-Baxter equation $R_1R_2R_1=R_2R_1R_2,$ where $R_1:=R\otimes I$, $R_2:=I\otimes R$, $I$ denotes the identity matrix of degree $d$.
\item[ii)] $R$ satisfies the Hecke equation $(R-q)(R+1)=0$ for some $q$ which will be assumed {\it not to be a root of unity} of order larger than 1.
\item[iii)] There exists a matrix $P_{ij}^{kl}$ such that $P_{jn}^{im}R^{nk}_{ml}=\delta^i_l\delta^k_j.$ \end{itemize}

\noindent{\sc Example.} The following main example of Hecke symmetries 
was first considered by Manin \cite{manin2}. Assume that the variables 
$x_i$, $i\leq r$ are even and the rest $s=d-r$ variables are odd. Define, for 
$1\leq i,j,k,l\leq r+s$,
$$ {R^{(r|s)}}^{kl}_{ij}:=\left\{
\begin{array}{lll}
q^2&\mbox{ if }& i=j=k=l, \hi=0\\
-1&\mbox{ if }& i=j=k=l, \hi=1\\
q^2-1 &\mbox{ if }&k=i<j=l\\
 (-1)^{\hat{i}\hat{j}}q&\mbox{ if }&k=j\neq i=l\\
 0&\multicolumn{2}{l}{\mbox{ otherwise. }}
\end{array}\right.$$
The Hecke equation for $R^{(r|s)}$ is $(x-q^2)(x+1)=0$.
When $q=1$, $R^{(r|s)}$ reduces to the super-permuting operator on $V\otimes V$.

  Let $\{z^i_j,t_j^i|1\leq i,j\leq d\}$ be a set of variables, where the parities of ${x^i_j}$ and ${t^i_j}$ are $\hi+\hj$. 

The super algebra $E_R$ to be the quotient algebra of the free non-commutative algebra on the generators $\{z^i_j|1\leq i,j\leq d\}$, by the relations
\begin{eqnarray}\label{eq1.1} (-1)^{\hs(\hi+\hp)}R^{kl}_{ps} z^p_iz^s_j&=&(-1)^{\hl(\hq+\hk)}z^k_qz^l_nR_{ij}^{qn},\quad 1\leq i,j,k,l\leq d.\end{eqnarray}
Here, we use the convention of summing up over the indices that appear in both lower and upper places.

The super algebra $H_R$ is defined to be the quotient of the free non-commutative algebra generated by $\{z^i_j,t_j^i|1\leq i,j\leq d\}$, by the relations
\begin{eqnarray}\label{eq1.2}(-1)^{\hs(\hi+\hp)}R^{kl}_{ps} z^p_iz^\sym_j&=&(-1)^{\hl(\hq+\hk)}z^k_qz^l_nR_{ij}^{qn},\quad 1\leq i,j,k,l\leq d,\\ 
\label{eq1.3}(-1)^{\hat{j}(\hj+\hk)}z^i_jt^j_k&=&(-1)^{\hat{l}(\hl+\hi)}t^i_lz^l_k=\delta^i_k,\quad 1\leq i,k\leq d.\end{eqnarray}

 The super algebra $E_R$ is a super bialgebra with the coproduct  given by 
$$\Delta(z^i_j)=z^i_k\otimes z^k_j,\quad \Delta(t^i_j)=t_j^k\otimes t_k^i.$$
 The super algebra $H_R$ is  a Hopf super algebra with the coproduct   given by
$$\Delta(z^i_j)=z^i_k\otimes z^k_j,\quad \Delta(t^i_j)=t_j^k\otimes t_k^i,$$
and  the antipode given by 
$$S(z^i_j)=(-1)^{\hj(\hi+\hj)}t^i_j,\quad S(t_j^i)=(-1)^{\hi(\hi+\hj)}C^i_kz^k_l{C^{-1}}^l_j,$$
 where $C^i_j:=P^{il}_{jl}$. See \cite{phhai01b} for details.
 
 The super bialgebra $E_R$ is called the (function algebra on) a quantum matrix super semigroup $\text{M}_q(r|s)$.
The Hopf  super algebra $H_R$ is called the (function algebra on) a quantum general linear group  GL$_q(r|s)$.
When
$R=R^{(r|s)}$ the associated Hopf super algebra is called the (function algebra on) standard quantum general linear super group GL$_q(r|s)$.
 Note that $R^{(r|s)}$ has birank $(r,s)$. 
 
The Hecke algebra of type $A$, $\H_n=\H_{n,q}$ is generated by elements $T_i,1\leq i\leq n-1$, subject to the relations
$$\begin{array}{l} T_iT_j=T_jT_i, |i-j|\geq 2;\\
T_iT_{i+1}T_i=T_{i+1}T_iT_{i+1};\\
 T_i^2=(q-1)T_i+q.\end{array}$$
To each element $w$ of the symmetric group $\mathfrak S_n$, one can associate in a canonical way an element $T_w$ of $\H_n$, in particular, $T_1=1, T_{(i,i+1)}=T_i$. The set $\{T_w|w\in\mathfrak S_n\}$ forms a $\mathbb k$ basis for $\H_n$.

The operator $R$ induces an action of the Hecke algebra $\mathcal H_n$ on the tensor powers $V^{\otimes n}$ of $V$,  $\rho_n(T_i)=R_i:=\id_V^{i-1}\otimes R\otimes\id_V^{n-i-1}.$ We shall therefore use the notation $R_w:=\rho(T_w)$. On the other hand,  $E_R$ coacts on  $V$ by $\delta(x_i)=x_j\otimes z^j_i$. 
Since $E_R$ is a bialgebra, it coacts on $V^{\otimes n}$ by means of its multiplication.  
With the assumption that $q$ is not a root of unity of order larger than 1, $\mathcal H_n$ is semi-simple and
we have the double centralizer theorem asserting that the action and coation mentioned here are  centralizers of each other in $\End_{\mathbb k}(V^{\otimes n})$ \cite{phhai01b}. It follows that $E_R$-comodules are semi-simple and each simple $E_R$-comodule is the image of the operator induced by a primitive idempotent of $\mathcal H_n$ and, conversely, each primitive idempotent of $\mathcal H_n$ induces an $E_R$ comodule which is either zero or simple. Since irreducible representations of $\mathcal H_n$ are parameterized by partitions of $n$, primitive idempotents of $\mathcal H_n$, up to conjugation, are parameterized by partitions of $n$, too.

For example, using the notation
$$[n]:=\frac{q^n-1}{q-1};\quad [n]!:=[1][2]\ldots[n],$$
we have the (central) primitive idempotents 
$$x_n:=\frac1{[n]!}\sum_w T_w\quad \mathrm{and}\quad  y_n:=\frac1{[n]!}q^{n(n-1)/2}\sum_w (-q)^{-l(w)}T_w,$$
which induce the symmetrizing and anti-symmetrizing operators $X_n$, resp. $Y_n$, on $V^{\otimes n}$. Let $\sym_n:=\Im X_n$ and $\exterior_n:=\Im Y_n$. 
One can show that $ \sym_n$ (resp. $\exterior_n$) is isomorphic to the $n$-th homogeneous compoment of the quadratic  algebra $\sym(V)$ (resp. $\exterior(V)$) defined as follows: 
$$\sym\cong T(V)/(\Im(R-q)),\quad (\text{resp. } \exterior \cong T(V)/(\Im(R+1))),$$
($T(V)$ denotes the tensor super algebra on $V$).
 These algebras are called the {\em symmetric and exterior tensor algebras} on a quantum super space. 

By definition, the Poincar\'e series $P_\exterior(t)$ of $\exterior$ is $\sum_{n=0}^\infty\dim_{\mathbb k}(\exterior_n) t^n$. It is proved that this series is a rational function having only real negative roots and real positive poles \cite{phhai99}. Let $r$ be the number of its roots and $s$ be the number of its poles. Then simple $E_R$-comodules are parameterized by hook-partitions from  $\Gamma^{rs}_n:=\{\lambda\vdash n|\lambda_{r+1}\leq s\}$~\cite{phhai01b}.

Simple $H_R$-comodules are much more complicated. The main difficulty lies in the fact that $H_R$-comodules are not semi-simple. In \cite{phhai05} it is shown that, as a braided monoidal category, the category of $H_R$-comodules depends only on the quantum parameter $q$ and the birank of $R$. Thus the problem reduces to the case of the standard deformation $R^{(r|s)}$. In this case the problem was studied by R.B. Zhang, et.al. \cite{zhang93,zhang06}, using the duality between $H_{R^{(r|s)}}$ and $\mathcal U_q(\frak{gl}(r|s))$. 

The problem of constructing all its simple comodules is still open.  The aim of this work is to treat this problem in the particular case, when $R$ has birank $(3,1)$.

 \section{The double Koszul complex}
  \subsection{The Koszul complex $K$}
  The Koszul complex $K$ associated to $R$ can be defined as  a collection of complexes $K_a$. 
  The terms of $K_a$ are indexed by pairs $(k,l)$ with $k-l=a$. Denote by $\db:\mathbb k\to V\otimes V^*$ the map
  $1\mapsto x_i\otimes \xi^i$, where $(\xi^i)$ is the basis of $V^*$, dual to the basis $(x_i)$ of $V$.
The  term $K_{k,l}$ is $\exterior_k\otimes \sym_{l}{}^*$ and 
  the differential $d_{k,l}:\exterior_k\otimes \sym_{l}{}^*\to \exterior_{k+1}\otimes \sym_{l+1}{}^*$ is given by: 
$$ d_{k,l}:\exterior_k\otimes \sym_{l}{}^*\longrightarrow V^{\otimes k}\otimes V^{*\otimes l}
\stackrel{\id\otimes\db\otimes \id}{\longrightarrow} 
V^{\otimes k+1}\otimes V^{*\otimes l+1}\stackrel{Y_{k+1}\otimes X_{l+1}{}^*}{\longrightarrow}
\exterior_{k+1}\otimes \sym_{l+1}{}^*,$$
where $X_l, Y_k$ are the $q$-symmetrizing operators introduced in 
Section \ref{sect1}.  The reader is referred to \cite{gur1} for the proof that
$d$ is a differential.  

Define the maps $\partial_{k,l}$ as follows:
$$ \partial_{k,l}:\exterior_{k+1}\otimes \sym_{l+1}{}^*\rightarrow
 V^{\otimes k+1}\otimes V^{*\otimes l+1}\stackrel{\scriptscriptstyle \id
\otimes(\ev R_{V, V^*})\otimes\id}{\longrightarrow} V^{\otimes k}\otimes V^{*\otimes l}
\stackrel{\scriptscriptstyle Y_{k}\otimes X_{l}{}^*}{\rightarrow}\exterior_{k}\otimes \sym_{l}{}^*,$$
where $\ev:V^*\otimes V\to \mathbb k$ is the evaluation map and $R_{V,V^*}:V\otimes V^*\to V^*\otimes V$ is the symmetry induced from $R$. In terms of the dual bases $(x_i)$ and $(\xi^j)$ it is given by $x_i\otimes \xi^j\mapsto \xi^k \otimes x_l P^{jl}_{ik}$, thus $\ev R_{V,V^*}(x_i\otimes \xi^j)=C^i_j$.

One can show \cite{gur1,phhai05} that $\partial$ is also a differential and satisfies
\begin{equation}\label{ct3}
q[l][k]d\partial+[l+1][k+1]\partial d=q^k([l-k]-[r-s])\id
\end{equation}
on ${K_{k,l}}$, where $(r,s)$ is the birank of $R$.
Consequently, the complex $K_a$ is exact if $a\neq s-r$. Further, it is shown that, for $a=s-r$, the complex $K_a$ is exact everywhere, except at the term $K_{r,s}$, where it has the one dimensional homology group.

\subsection{The Koszul Complex $L$} There is another Koszul complex associated to $V$, which
was first defined by Priddy as a free resolution of the symmetric tensor
algebra of $V$ (see \cite{Manin}). As in the case of the complex $K$, the complex $L$ is a collection of complexes $L_a$. The complex $L_a$ has $(p,r)$-term, with $p+r=a$, $L_{p,r} := \sym_p\otimes \exterior_r$ 
and  differential $ P_{p,r}: L_{p,r} \longrightarrow L_{p-1,r+1}$ 	
 given by
$$ \xymatrix @C=4em{P_{p,r}: \sym_p\otimes  \exterior_r \ar@{^(->}[r] &V^{\otimes p}\otimes  V^{\otimes r}  
\ar[r]^{  X_{p-1} \otimes Y_{r+1}}&\sym_{p-1}\otimes  \exterior_{r+1}}.$$
The complexes $(L_a , P), a\geq 1,$ are exact.
This is shown by considering the map $Q_{p,r}: L_{p-1,r+1}\longrightarrow L_{p,r} $,   given by
$$ \xymatrix@C=4em { Q_{p,r}: \sym_{p-1}\otimes \exterior_{r+1}\ar @{^(->}[r]&V^{\otimes p-1}\otimes V^{\otimes r+1}   \ar[r]^{X_p \otimes Y_r}& \sym_p\otimes \exterior_r}.$$
 One checks \cite{gur1} that on $ L_{p,r}$ 
 \begin{equation}\label{ct60}
[r][p+1] PQ + [p][r+1]QP =[p+r]\id.
\end{equation}   
\begin{rmk}\em
 The differentials of both complexes are morphisms of $H_R$-comodules.
\end{rmk}
\subsection{The double Koszul complex}
The two Koszul complexes mentioned in the previous section can be combined into a double complex called the double Koszul complex.  
For simplicity we shall  use the dot ``$ \cdot $'' to denote the tensor product.  Fix  an integer $a$.
 We arrange the Koszul complexes $K_{-a}, K_{-a-1}, K_{-a-2},$...  as follows. 
 $$ \xymatrix@C=4ex@R=2ex{ 
K_{-a}\ :0\ar[r]&S_{a}{}^*\ar[r]^{d_{0,a}}&\exterior_1\cdot \sym_{a+1}{}^*\ar[r]^{d_{1,a+1}}&\exterior _2\cdot \sym_{a+2}{}^*\ar[r]^{d_{2,a+2}}
&\exterior _3\cdot \sym_{a+3}{}^*\ar[r]&\ldots\\ 
K_{-a-1}: \ &0\ar[r]&\sym_{a+1}{}^*\ar[r]^{d_{0,a+1}}&\exterior_1\cdot \sym_{a+2}{}^* \ar[r]^{d_{1,a+2}}&\exterior _2\cdot \sym_{a+3}{}^*\ar[r]&\ldots\\ 
K_{-a -2}:\  &&0\ar[r]&\sym_{a+2}{}^*\ar[r]^{d_{0,a+2}}&\exterior_1\cdot \sym_{a+3}{}^* \ar[r]&\ldots }$$
here $\sym_i$ and $\exterior_i$ are set to 0 if $i<0$.
To get the entries on a column into a complex we tensor  each complex $K_{i}$ with $\sym_{-a-i}$,
i.e. the complex $K_{-1-a}$ is tensored with $\sym_1$, the complex $K_{-2-a}$ is tensor with $\sym_2$...
Then each  column can be interpreted as the complexes
$L_j$ tensored with $\sym_{a+j}{}^*$.  
 Thus we have the following diagram with all rows being the Koszul complexes $K_\bullet$ tensored with
$\sym_\bullet$ and columns are the Koszul complexes $L_\bullet$ tensored with $\sym_\bullet{}^*$:
\begin{equation}\label{phuckep0}
\xymatrix@R=1.5em{& 0&0&0&0\\
0\ar@<.5ex>[r]    &\ar[u]\sym_{a}{}^* \ar@<.5ex>[r]^d   
&\ar[u]\exterior_1\cdot \sym_{a+1}{}^* \ar@<.5ex>[r]^d     &\ar[u]\exterior _2\cdot \sym_{a+2}{}^* 
\ar@<.5ex>[r]^d  
  &\ar[u]\exterior _3\cdot \sym_{a+3}{}^* \ar@<.5ex>[r]^d      &
 \ldots\\
   &0\ar@<.5ex>[r]  \ar@<1ex>[u] &\sym_1\cdot \sym_{a+1}{}^* \ar@<.5ex>[r]^d  
\ar@<1ex>[u]^{ P}  
&\sym_1\cdot \exterior _1\cdot \sym_{a+2}{}^* \ar@<.5ex>[r]^d   \ar@<1ex>[u]^{ P} &\sym_1\cdot \exterior _2\cdot \sym_{a+3}{}^*
\ar@<.5ex>[r]^d  
\ar@<1ex>[u]^{  P  }
&
  \ldots \\
 &&0\ar@<.5ex>[r]  \ar@<1ex>[u] &
  \sym_2\cdot \sym_{a+2}{}^* \ar@<.5ex>[r]^d    \ar@<1ex>[u]^{  P} 
&\sym_2\cdot \exterior _1\cdot \sym_{a+3}{}^*\ar@<.5ex>[r]^d  
\ar@<1ex>[u]^{  P}   & 
  \ldots \\
&&&0\ar[u]&\vdots\ar[u]}\end{equation}
  A general square in  diagram (\ref{phuckep0}) has the form
\begin{equation}\label{cell}\xymatrix@R=1.5em{
\sym_i\cdot \exterior_k\cdot \sym_l{}^*\ar[r]^{\id\otimes d}& \sym_i\cdot \exterior_{k+1}\cdot S_{l+1}{}^*\\
\sym_{i+1}\cdot \exterior_{k-1}\cdot S_l{}^*\ar[r]_{\id\otimes d}\ar[u]^{ P\otimes\id}&
\sym_{i+1}\cdot \exterior_k\cdot \sym_{l+1}{}^*\ar[u]_{  P\otimes\id} &\\
}\end{equation} 
with $l = i+k+a$.
For convenient, we   denote $ d:= \Id \otimes d,  P : =  P \otimes \Id$.  
It is easy to show that  $  P d = d  P$ for all these squares. Thus \eqref{phuckep0} is a bicomplex.

  We also  have an exact double Koszul complex with $d,  P$ replaced by $\partial,Q$.  
 {\small
 \begin{equation}\label{phuckep1}
\xymatrix@R=1.5em{& 0 \ar@{ -->}[d]  &0\ar@{ -->}[d]&0\ar@{-->}[d]&0\ar@{-->}[d]&
\\
 0  \ar@<-.5ex>@{<--}[r]_{\partial} &\sym_{a}{}^* \ar@{ -->}[d]^Q   \ar@<-.5ex>@{<--}[r]_{\partial}
& \exterior_1\cdot \sym_{a+1}{}^* \ar@{ -->}[d]^Q\ar@<-.5ex>@{<--}[r]_{\partial}     & \exterior _2\cdot \sym_{a+2}{}^* \ar@{ -->}[d]^Q
\ar@<-.5ex>@{<--}[r]_{\partial}
  & \exterior _3\cdot \sym_{a+3}{}^*\ar@<-.5ex>@{<--}[r]_{\partial}  \ar@{ -->}[d]^Q    & 
 \ldots\\
   &0\ar@<-.5ex>@{<--}[r]_{\partial}  &\sym_1\cdot \sym_{a+1}{}^* \ar@<-.5ex>@{<--}[r]_{\partial}
\ar@{ -->}[d]^Q
&\sym_1\cdot \exterior _1\cdot \sym_{a+2}{}^*\ar@<-.5ex>@{<--}[r]_{\partial}  \ar@{ -->}[d]^Q &\sym_1\cdot \exterior _2\cdot \sym_{a+3}{}^*
\ar@<-.5ex>@{<--}[r]_{\partial} 
\ar@{ -->}[d]^Q
& 
  \ldots \\
 &&0\ar@<-.5ex>@{<--}[r]_{\partial} &\ar@{ -->}[d]
  \sym_2\cdot \sym_{a+2}{}^*\ar@<-.5ex>@{<--}[r]_{\partial}    
&\ar@{ -->}[d] \sym_2\cdot \exterior _1\cdot \sym_{a+3}{}^*\ar@<-.5ex>@{<--}[r]_{\partial} 
 &  
  \ldots\\ &&&0 &\vdots }
  \end{equation}}
  
 
 From now, we assume that $R$ has birank $(3|1)$.

We combine the two diagrams \eqref{phuckep0} and
\eqref{phuckep1} into one: 
 
   \begin{equation}\label{phuckep3}
   \xymatrix{\sym_{i-1}\cdot \sym_{a+i-1}{}^*\ar@<.5ex>@{^(->}[r]^{d_{0,a+i-1}}   \ar@<-.5ex>@{<<--}[r]_{\partial_{0,a+i-1}}&
\sym_{i-1}\cdot \exterior_1\cdot \sym_{a+i}{}^* \ar@{-->>}[d]^Q\ar@<.5ex>[r]^{d_{1,{a+i}}}   \ar@<-.5ex>@{<--} [r]_{\partial_{1,{a+i}}}
&\sym_{i-1}\cdot \exterior_2 \cdot \sym_{a+i+1}{}^*\ar@{-->}[d]^Q\ar@<.5ex> [r]^{d_{2,{a+i}+1}}   \ar@<-.5ex>@{<--} [r]_{\partial_{2,{a+i}+1}}&  \sym_{i-1}\cdot \exterior_3 \cdot \sym_{a+i+2}{}^*\ar@{-->}[d]^Q
\cdots\\
          &\sym_i.\sym_{a+i}{}^* \ar@<1ex>@{^(->}[u]^{  P}\ar@<.5ex>@{^(->}[r]^{d_{0,{a+i}}}   \ar@<-.5ex>@{<<--}[r]_{\partial_{0,{a+i}}}& 
 \sym_i.\exterior _1.\sym_{{a+i}+1}{}^*  
\ar@<1ex>[u]^{  P}\ar@<.5ex>[r]^{d_{1,{a+i}+1}}  \ar@<-.5ex>@{<--}[r]_{\partial_{1,{a+i}+1}}   \ar@{-->>}[d]^Q&
\sym_i\cdot \exterior _2\cdot \sym_{{a+i}+2}{}^*   \ar@{-->}[d]^Q\ar@<1ex>[u]^{  P}\cdots\\
           & & \sym_{i+1}.\sym_{{a+i}+1}{}^* \ar@<1ex>@{^(->}[u]^{  P} 
\ar@<.5ex>@{^(->}[r]^{d_{0,{a+i}+1}}  \ar@<-.5ex>@{<<--}[r]_{\partial_{0,{a+i}+1}}\ar@<1ex>[u]^{ P}
& \sym_{i+1}\cdot \exterior _1\cdot \sym_{{a+i}+2}{}^*  \ar@<1ex>[u]^{  P}\cdots}
\end{equation}   

\begin{pro}\label{mde1}Assume that the Hecke symmetry $R$ has birank $(3,1)$.
Then the composed map   $  \partial    PQd  : \sym_i\cdot \sym_{{a+i}}{}^* \longrightarrow \sym_i\cdot \sym_{{a+i}}{}^* $ in   diagram  (\ref{phuckep3}) is an  isomorphism for all
 $a,i$ with $i,a+i\geq 0 $. Consequently $\sym_i\cdot \sym_{a}{}^*$ is isomorphic to a direct summand of $\sym_{i+1}\cdot \sym_{a+1}{}^*$. Moreover this isomorphism is an isomorphism of $H_R$-comodules.
   \end{pro} 
  {\it Proof. } We will use induction on $i$ to prove that the endomorphism $\partial   PQd : \sym_{i}\cdot \sym_{{a+i}}{}^* \longrightarrow  \sym_{i}\cdot \sym_{{a+i}}{}^*$ 
is diagonalizable with the set of  eigenvalues equal to
\begin{equation}\label{ctmd1}
 A_i :=\Big \{\frac{([a+2i+1-j]-[-2])[j]}{[i+1][a+i+1]}, j=1,2,\ldots,i+1
\Big \}
\end{equation} 
For  $i=0$, the map $PQ: \sym_a{}^* \longrightarrow \sym_a{}^*$ is equal  to $\id_{\sym_a{}^*}$. Hence 
$$g=\partial PQd= \frac {[a]-[-2]}{[a+1]}\id.$$

 Assume that  the claim holds true for  $i-1$. We have 
 \begin{align*}
  h:= \partial PQd  = & \partial [\frac {[i+1]-[2][i]QP}{[i+1]}]d = \partial d - \frac{[2][i]}{[i+1]}\partial QPd\\
 =& \partial d - \frac {[2][i]}{[i+1]}Q\frac {[q([a+i-1]-[-2])-[a+i]d \partial]}{[2][a+i+1]}P\\
 =& \partial d - \frac{q[i][([a+i-1]-[-2])]}{[i+1][a+i+1]}QP+ \frac{[i][a+i]}{[i+1][a+i+1]}Qd \partial P\\
 =& \Big[ \frac{[a+i]-[-2]}{[a+i+1]}-\frac{q[i]([a+i-1]-[-2]))}{[i+1][a+i+1]} \Big] \id + \frac{[i][a+i]}{[i+1][a+i+1]}Qd \partial P.
 \end{align*}
  By assumption
$  \partial PQd : \sym_{i-1}\cdot \sym_{a-1}{}^*\longrightarrow \sym_{i-1}\cdot \sym_{a-1}{}^* $ is diagonalizable with  eigenvalues in $A_{i-1},$, in particular it is invertible. Thus the minimal polynomial $P(X)$ of this operator has no multiple roots. It follows that the minimal polynomial of the operator
  $Qd\partial P: \sym_i\cdot \sym_{a+i}{}^*\longrightarrow  \sym_i\cdot \sym_{a+i}{}^*$ is just $XP(X)$.
Consequently $Qd\partial P$ is diagonalizable with   eigenvalues  in $A_{i-1}\cup \{ 0\}.$
Thus $ \partial   PQd: \sym_{i}\cdot \sym_{a}{}^*\longrightarrow \sym_{i}\cdot \sym_{a}{}^*$ 
is diagonalizable with the set   eigenvalues in $A_i$ .$\Box$

 Consider the diagram in (\ref{phuckep0}) as a exact sequence of horizontal complexes (except for the first column) and split it into short exact sequences.

\begin{equation}\label{phuckep4}
  \xymatrix{
 \ldots \Ker P_{i,k}\cdot S_{i+k+a}{}^*\ar@<.5ex>[r]^{d'_{k,i+k+a}}\ar@{-->}[d]^Q&  \Ker P_{i,k+1}\cdot S_{i+k+a+1}{}^*\ar@{-->}[d]^Q
\ar@<.5ex>[r]^{d'_{k+1,i+k+a+1}}&  \Ker P_{i,k+2}\cdot S_{i+k+a+2}{}^*\ar@{-->}[d]^Q
\ldots \\
\ldots\sym_{i+1}\cdot \exterior _{k-1} \cdot S_{i+k+a}{}^*\ar@{-->}[d]^Q\ar@<1ex>@{->>}[u]^{  P_{i+1,k-1}} 
  \ar@<.5ex>[r]^{d_{k-1,i+k+a}}
 &\sym_{i+1}\cdot \exterior _k\cdot S_{i+k+a+1}{}^*\ar@{-->}[d]^Q
\ar@<1ex>@{->>}[u]^{ P_{i+1,k}}   \ar@<.5ex>[r]^{d_{k,i+k+a+1}}
 &\sym_{i+1}\cdot \exterior _{k+1}\cdot S_{i+k+a+2}{}^*\ar@<1ex>@{->>}[u]^{  P_{i+1,k+1}}\ar@{-->}[d]^Q\ldots\\
\ldots\Ker P_{i+1,k-1}\cdot S_{i+k+a}{}^* \ar@<.5ex>[r]^{d'_{k-1,i+k+a}} \ar@<1ex> @{^(->}[u]^{i}&\Ker P_{i+1,k }
\cdot \sym_{i+k+a+1}{}^* 
\ar@<.5ex>[r]^{d'_{k,i+k+a+1}}  \ar@<1ex> @{^(->}[u]^{i}&
\Ker P_{i+1,k+1}\cdot S_{i+k+a+2}{}^*\ar@<1ex>@{^(->}[u]^{i}\ldots }
\end{equation}
where $d'_{k,i+k+a}$ is the restriction of $d_{k,i+k+a}$ to $\Ker P_{i,k}\cdot \sym_{i+k+a}{}^*$. Notice that $\Ker P_{i,j} = \Im P_{i+1,j-1}$ for all $i\geq 0.$


    Consider the following part of \eqref{phuckep4}  for $i, k\geq 1$:
{\small
\begin{equation}\label{sodo}\xymatrix{ 
& 
\sym_i\cdot \exterior_{k+1} \cdot \sym_{a+i+k+1}{}^*\ar[r] \ar @<1 ex>@{-->}[d]& \cdot \sym_i \cdot \exterior_{k+2}\cdot S_{a+i+k+2}{}^*
\ar@<1ex>@{-->}[l] \ar @<1 ex>@{-->}[d] \\
&  \Ker P_{i,k+1}\cdot S_{a+i+k+1}{}^*\ar[r] \ar @{^(->}[u] \ar @<1ex>@{-->}[d]^Q&
\Ker P_{i,k+2}\cdot \sym_{a+i+k+2}{}^*\ar @{^(->}[u] \ar @<1ex>@{-->}[d] \ar@<1ex>@{-->}[l]   \\
&  \sym_{i+1}\cdot \exterior_{k}\cdot S_{a+i+k+1}{}^*\ar[r]^d  \ar  @{->>}[u]^{ P}&
\sym_{i+1}\cdot \exterior_{k+1}\cdot \sym_{a+i+k+2}{}^*  \ar @{->>}[u] \ar@<1ex>@{-->}[l]^{\partial} }
\end{equation}}

 \begin{pro}\label{dl1}
Assume that the Hecke symmetry $R$ has birank $(3,1)$. 
Then for $i\geq 0, k\geq 1,a+i+k+1\geq 0$ the composed map
 $$ P\partial dQ: \Ker P_{i,k+1}\cdot \sym_{a+i +k+1}{}^*\longrightarrow \Ker P_{i,k+1}\cdot \sym_{a+i+k+1}{}^*$$ 
 in  diagram  (\ref{sodo}) 
 is an  isomorphism.  Consequently $\Ker P_{i,k+1}\cdot \sym_{a+i+k+1}{}^*$ is isomorphic to a direct summand of $\sym_{i+1}\cdot \Im d_{k,a+i+k+1}.$ Moreover the isomorphism is an isomorphism of  $H_R$-comodules.
\end{pro}

{\it Proof.} We assume first that $a\geq 0$, the case $a<0$ is treated similarly but a bit more tedious.
 We use induction to prove that 
$$ P\partial dQ: \Ker P_{i,k+1}\cdot \sym_{a+i +k+1}{}^*\longrightarrow \Ker P_{i,k+1}\cdot \sym_{a+i+k+1}{}^*$$
is diagonalizable with  eigenvalues  
 
\begin{align*}
&A_i := \Big \{\frac{q^k([a+k+2i-j+2]-[-2])[j]}{[i+1][k+1]^2[a+i+k+2]}, j=1,2,\ldots,i+1,i+k+1
 \Big \}
\end{align*}
 
 For $i =0$, consider  the following part of \eqref{sodo}:
 {\small
\begin{equation*} \xymatrix{\exterior_{k }\cdot   \sym_{a+k }{}^*\ar[r]_d \ar @<1 ex>@{-->}[d]^Q& \ar@<-.5ex>@{-->}[l]_{\partial}
  \exterior_{k+1}\cdot   \sym_{a+k+1}{}^*\ar[r]_d \ar@<1ex>@{-->}[d]^Q& 
\cdot \exterior_{k+2}\cdot S_{a+k+2}{}^*\ar@<1ex>@{-->}[d]^Q \ar@<-.5ex>@{-->}[l]_{\partial} \\
  \sym_1\cdot\exterior_{k-1}\cdot S_{a+k}{}^*\ar[r]_d \ar@{->>}[u]^{  P} & \sym_{1}\cdot \exterior_{k}\cdot S_{a+k+1}{}^*\ar[r]_d 
\ar@{->>}[u]^{  P}
\ar@<-.5ex>@{-->}[l]_{\partial}&
\sym_{1}\cdot \exterior_{k+1}\cdot \sym_{a+k+2}{}^*\ar@{->>}[u]^P   \ar@<-.5ex>@{-->}[l]_{\partial} 
}\end{equation*} }
The composed map $ P\partial DQ: \exterior_{k+1}\cdot \sym_{a+k+1}{}^*\longrightarrow \exterior_{k+1}\cdot \sym_{a+k+1}{}^*.$
By means of   formulas (\ref{ct3}) and (\ref{ct60}) we have 
\begin{eqnarray*}
P\partial dQ &=& P\frac{[q^k([a+1]-[-2])- [k][a+k+1]d\partial ]}{[k+1][a+k+2]} Q  \\ &= &\frac {q^k([a+1]-[-2])} {[k+1][a+k+2]}\Id   - \frac {[k][a+k+1]}{[k+1][a+k+2]}d \partial\\
\end{eqnarray*}
Since  $d\partial $  is diagonalizable with eigenvalues    $0$ and $ \frac {[a]-[-2]}{ [k+1][a+k+1]}$,  
  $ P\partial dQ$ is diagonalizable with the set of   eigenvalues    
$$ A_0 := \Big \{\frac {q^k[k+1]([a+1]-[-2])}{[k+1]^2[a+k+2]}, 
 \frac {q^k([a+k+1]-[-2])}{[k+1]^2[a+k+2]}\Big \}.$$
 For $i=1$, 
 consider diagram \eqref{sodo} with $i=1$ the map $ P\partial dQ: \Ker P_{1,k+1}\cdot \sym_{a+k+2}{}^*\longrightarrow \Ker P_{1,k+1}\cdot \sym_{a+k+2}{}^*$, we have

\begin{eqnarray*}
 P\partial dQ &=& P\frac{[q^k([a+2]-[-2])-q[k][a+k+2]d \partial]}{[k+1][a+k+3]}Q \\
&=& \frac {q^k([a+2]-[-2])}{[k+1][a+k+3]} PQ -\frac{q[k][a+k+2]}{[k+1][a+k+3]}  dPQ\partial \\
&=&\frac{q^k([a+2]-[-2])[k+2]}{[2][k+1]^2[a+k+3]}id - \frac{q[k][a+k+2]}{[k+1][a+k+3]}d[\frac{[k+1]-[k+1]QP}{[2][k]}] \partial \\
&=& \frac{q^k([a+2]-[-2])[k+2]}{[2][k+1]^2[a+k+3]}id - \frac{q[a+k+2]}{[2][a+k+3]}d \partial + \frac{q[a+k+2]}{[2][a+k+3]}dQP\partial
 \end{eqnarray*} 
We have $d \partial: \sym_1 \cdot \exterior \cdot \sym_{a+k+2}{}^*\longrightarrow  \sym_1 \cdot \exterior \cdot \sym_{a+k+2}{}^*$ is diagonalizable with eigenvalues 
$$\frac{q^{k+1}([a+1]-[-2])}{q[k+1][a+k+2]}\quad \text{ and } 0.$$
 On the other hand, we have $d \partial \cdot dQP \partial = dQP \partial \cdot d \partial$
and if $d \partial(x)=0$ than  $dQP \partial(x)=0.$ 
Therefore, the eigenvalues of   $P\partial dQ : \Ker P_{1,k+1}\cdot \sym_{a+k+2}{}^* \longrightarrow  \Ker P_{1,k+1}\cdot \sym_{a+k+2}{}^* $ are in the set
 \begin{align*}
A_1 := &  \Big \{\frac {q^k([a+2]-[-2])[k+2]}{[2][k+1]^2[a+k+3]}, \frac{q^k([a+k+3]-[-2])}{[2][ k+1]^2[a+k+3]}, \frac {q^k[2]([a+k+2-[-2])}{[2][k+1]^2[a+k+3]}
 \Big \}. 
\end{align*}

  In general, consider the composed map $$P\partial dQ: \Ker P_{i,k+1}\cdot \sym_{a+i+k+1}{}^* \longrightarrow \Ker P_{i,k+1}\cdot \sym_{a+i+k+1}{}^*
$$ in diagram (\ref{sodo}), we have
 
\begin{align*}
P\partial dQ =& P\left[\frac{q^k([a+i+1]-[-2]) - q[k][a+i+k+1]d\partial }{[k+1][a+i+k+2]}\right]Q\\
=& \frac {q^k([a+i+1]-[-2])PQ}{[k+1][a+i+k+2]}- \frac{q[k][a+i+k+1]dPQ}{[k+1][a+i+k+2]}\partial \\
=& \frac{q^k([a+i+1]-[-2])[i+k+1]\Id}{[k+1]^2[i+1][a+i+k+2]}\\
&\quad -\frac{q[k][a+i+k+1]d}{[k+1][a+i+k+2]}\cdot\frac {([i+k] - [i][k+1]QP)\partial}{[k][i+1]} \\
=&  \frac{q^k([a+i+1]-[-2])[i+k+1]\Id }{[k+1]^2[i+1][a+i+k+2]} - \frac{q[i+k][a+i+k+1]d\partial}{[k+1][i+1][a+i+k+2]} \\
&\quad + \frac{q[i][a+i+k+1]dQP\partial }{[i+1][a+i+k+2]}.
\end{align*}
One has $d \partial$ is diagonalizable with the set of eigenvalues 
$$\left\{\frac{q^k([a+i]-[-2])}{[k+1][a+i+k+1]}, 0\right\}.$$ 
We have $d \partial \circ dQP \partial = d QP \partial \circ d \partial$ and if $d \partial(x)=0 $ then $dQP \partial (x)=0.$ By induction assumption
  $P\partial d Q: \Ker P_{i-1,k+1}\cdot \sym_{a+i+k}{}^* \longrightarrow 
 \Ker P_{i-1,k+1}\cdot \sym_{a+i+k}{}^*$  is diagonalizable with eigenvalues in the set $A_{i-1}$. Thus  
the composed map  
$$  dQP \partial: \Ker P_{i,k+1}\cdot \sym_{a+i+k+1}{}^* \longrightarrow 
 \Ker P_{i,k+1}\cdot \sym_{a+i+k+1}{}^*$$
  is diagonalizable with  the set of   eigenvalues is  $ A_i$. The proof is complete.$\Box$

 \section{Construction of irreducible representations of $GL_q(V)$.}
Let $R: V \otimes V \rightarrow V \otimes V$ be a Hecke symmetry with birank $(3|1)$. 
 Using the double Koszul complex, we will construct in this section for each (integrable) dominant weight, i.e. a quadruple $(m,n,p,t)$ of integers, with $m\geq n\geq p$, a comodule $I(m,n,p|t)$ of $H_R$.  The proof that these comodules are simple and furnish all simple $H_R$-comodules will be given in the next section.
 
Recall that the complex $K_2$ is exact  everywhere, except at the  term $K_{3,1}$, where the homology is one dimensional. Denote this comodule  by $I(1,1,1|1)$.

 For a dominant weight $(m,n,p|t)$ set 
$$I(m,n,p|t):= I(m-t,n-t,p-t|0) \otimes I(1,1,1|1){}^{\otimes t}.$$
Thus one is led to construct $I(m,n,p|0)$.

First, recall from Section \ref{sect1} that each partition $\lambda\in \Gamma^{3|1}$ defines a simple $H_R$-comodules. Denote it by $M_\lambda$. Such a partition $\lambda$ has the form $(\lambda_1,\lambda_2,\lambda_3,1^{\lambda_4})$. For a weight $(m,n,p|0)$ with $p\geq 0$ set
\begin{equation}\label{I1}
I(m,n,p|0):=M_{(m,n,p)}.
\end{equation}
Further, for such a dominant weight with $p\geq 1$ we set
\begin{equation}\label{I2}\begin{array}{ll}
I(-p-2,-n-2,-m-2|0):=I(m,n,p|0){}^*\otimes \ I(1,1,1|1){}^{*\otimes 3}\end{array}
\end{equation}
and for a dominant weight of  type $(m,n,0|0)$ we set
\begin{equation}\label{I3}
\begin{array}{ll}
 I(-2,-n-1,-m-1|0):=I(m,n,0|0){}^*\otimes \ I(1,1,1|1){}^{*\otimes 2}.\end{array}
\end{equation}
Finally we set 
\begin{equation}\label{I4}\begin{array}{ll}
I(-1,-1,-m|0):=I(m,0,0|0){}^*\otimes\ I(1,1,1|1){}^{*\otimes 1}\end{array}
\end{equation}
The reason for the choice of the weight on the left hand side above will be explained in the next section when we compute the character.
 
   \subsection{Comodules constructed from complex $K$}
Consider complexes  $K_a$, with  $a: = k-l  \neq 2$.
$$K_a: \ldots\longrightarrow\exterior _{k-1}\otimes\sym_{l-1}{}^*\longrightarrow \exterior _k\otimes\sym_l {}^*\longrightarrow  \exterior _{k+1}\otimes\sym_{l+1}{}^*\longrightarrow  \ldots$$
 By using the exactness of the complex $K$ we will construct a class of irreducible representations of $GL_q(3|1)$. According to (\ref{ct3}) we have  
\begin{equation}\label{xdanh}
 \exterior _k.\sym_l{}^* \cong  \Im d_{k-1,l-1} \oplus \Im d_{k,l}.
 \end{equation}


For a dominant weight $(m,m,p|0)$ with $m\geq 0>p$, set  
\begin{equation}\label{I5}
I(m,m,p|0) := \Im d_{m+2,m-p} \otimes I(1,1,1|1){}^{\otimes m-1},
\end{equation} 
and
\begin{equation}\label{I6}
\begin{array}{l}
I(-2-p,-m-2,-m-2|0):=I(m,m,p|0){}^*\otimes I(1,1,1|1){}^{*\otimes 3}.
\end{array}
\end{equation}

  \subsection{Comodules constructed from the double Koszul complex.}
From Proposition~\ref{mde1}, for any $i,a$ with $i,a+i\geq 0$, there exists $X_{i,a}$ such that 
$$\sym_{i+1}. \sym_{a+i+1}{}^*= \sym_{i}.\sym_{a+i}{}^* \oplus X_{i,a}.$$
For any dominant weight $(m,-1,p|0)$  with $m\geq 0$ (and $p\leq-1$), set
\begin{equation}\label{I7}
I(m,-1,p|0)= X_{m,-m-p-1}\otimes I(1,1,1|1){}^*.
\end{equation}

  According to Proposition \ref{dl1}, there exists comodule $Y_{i,k,a}$ such that, for $i,k,a$ with $k\geq 1$, $i,a+i+k+1\geq 0$,
  $$\Ker P_{i,k+1}\otimes S_{a+i+k+1}{}^*\oplus Y_{i,k,a}\cong
  S_{i+1}\otimes \Im d_{k,a+i+k+1}.$$
    For a dominant weight $(m,n,p|0)$ with $m>n\geq 0>p$, set
  \begin{equation}\label{I8}
I(m,n,p|0)=Y_{m-n-1,n+2,n-m-p-2}\otimes I(1,1,1|1)^{*\otimes n-1}.\end{equation}
For a dominant weight $(m,n,p|0)$ with $m\neq -2, n\leq -2$,
we set
\begin{equation}\label{I9}
I(m,n,p|0)=I(-2-p,-2-n,-2-m){}^*I(1,1,1|1)^{*\otimes 3}.
\end{equation}

Thus for any integrable dominant  weight $(m,n,p|0)$ we have constructed a comodule $I(m,n,p|0)$. Here is the detailed check:
\begin{enumerate}
\item  $m\geq n\geq 0$: $I(m,n,p|0)$ is given by \eqref{I1}.
\item  $m\geq n\geq 0> p$:
	\begin{enumerate}
	\item  $m=n$: $I(m,m,p|0)$ is given by \eqref{I5};
	\item  $m>n$: $I(m,n,p|0)$ is given by \eqref{I8};
	\end{enumerate}
\item  $m\geq 0>n\geq p$:
	\begin{enumerate}
	\item  $n=-1$: $I(m,-1,p|0)$ is given by \eqref{I7};
	\item  $-2\geq n$: $I(m,n,p|0)$ is given by \eqref{I9};
	\end{enumerate}
\item  $0>m\geq n\geq p$: 
	\begin{enumerate}
	\item  $m=n=-1$: $I(-1,-1,p|0)$ is given by \eqref{I4};
	\item  $m=-1>n$: $I(-1,n,p|0)$ is given by \eqref{I9};
	\item  $m=-2$: $I(-2,n,p|0)$ is given by \eqref{I3};
	\item  $-2>m$: $I(m,n,p|0)$ is given by \eqref{I9}.
	\end{enumerate}
\end{enumerate}

In the next section we shall exhibit the simplicity of these comodules by reducing it to the case of the standard Hecke symmetry $R^{(r|s)}$ and using the formal character.

\section{Simplicity and completeness}
In this section we shall prove the simplicity of the comodules constructed in the previous section and that they furnish all simple comodules of $H_R$.
Our method is to use the representation theory of the quantum universal enveloping algebra ${\mathcal 
U}_q(\mathfrak{gl}(3|1))$. According to \cite[Thm~4.3]{phhai05} there is a monoidal equivalence between the category of comodules over $H_R$ and the category of comodules over $H_{R^{(3|1)}}$. Thus the problem is reduced to the case $R=R^{(3|1)}$. In this case, there is a duality between $H_{R^{(3|1)}}$ and ${\mathcal U}_q(\mathfrak{gl}(3|1))$, \cite[Thm~3.5]{zhang06}, which shows that there is an equivalence between the category of comodules over $H_{R^{(3|1)}}$ and finite dimensional integrable representations of ${\mathcal U}_q(\mathfrak{gl}(3|1))$. Notice that irreducible representations of $\mathcal U_q(\mathfrak{gl}(n|1))$ can be obtained by other methods, see e.g. \cite{palev, Ky}. But these methods are not compatible with the braided monoidal equivalence mentioned here. This is the reason why we want to give a construction based merely on the braiding (given by $R$) and the two maps $\ev$ and $\db$.

For finite dimensional representations of ${\mathcal U}_q(\mathfrak{gl}(3|1))$ the weight decomposition is obtained in the same way as for the classical case of $\mathfrak{gl}(3|1)$, whence the character is defined and does not depend on the parameter $q$ (as long as $q$ is not a root of unity).

The character of $H_{R^{(3|1)}}$-comodules can be defined directly. Consider the the quotient Hopf super algebra of this algebra by setting $z^i_j=0$ for all $i\neq j$. This quotient is just the algebra of Lorenz polynomials $\mathbb k[z^i_i{}^{\pm1}]$. Any $H_{R^{(3|1)}}$-comodule restricts to a comodule over $\mathbb k[z^i_i{}^{\pm1}]$, which is semi-simple, yielding the weight decomposition, which is independent of the quantum parameter. Assume that $M$ is a comodules over $H_{R^{(3|1)}}$, consider it as a comodule over  $\mathbb k[z^i_i{}^{\pm1}]$ we obtain the decomposition
$$M\cong \bigoplus_\lambda M_\lambda,$$
where $\lambda$ run over the set of $\mathbb Z$-linear mappings from the free abelian group generated by $z^i_i$ to $\mathbb Z$, i.e. the set of integrable weights. The character of $M_\lambda$ is defined to be
$$\ch(M_\lambda):=\sum \dim_{\mathbb k}(M_\lambda) e^\lambda.$$
If follows immediately from the definition that the character is additive with respect to short exact sequences and multiplicative with respect to the tensor product.
The fact that this definition agrees with the above definition follows from the explicit duality between $H_{R^{(3|1)}}$ and 
${\mathcal U}_q(\mathfrak{gl}(3|1))$.

Now to finish the proof that all comodules of $H_R$ constructed in the previous section are simple and furnish all $H_R$-comodules, it suffices to verify the following lemma and to compute explicitly the character of these comodules.
\begin{lem}\label{lem_weight}
Let $V$ be a representation of ${\mathcal U}_q(\mathfrak{gl}(3|1))$ with the character equal to the character of the simple highest weight representation $V(\lambda)$. Then $V$ is isomorphic to $V(\lambda)$.
\end{lem}
\begin{proof}
Consider $V$ and $V(\lambda)$ as representations of the Hopf subalgebra ${\mathcal U}_q(\mathfrak{gl}(3)\oplus \mathfrak{gl}(1))$. Since they have the same character, they are isomorphic. In particular, as ${\mathcal U}_q(\mathfrak{gl}(3)\oplus \mathfrak{gl}(1))$-representations, $V$ contains a direct summand with highest weight $\lambda$, say $S(\lambda)$.

According to \cite{zhang93}, $V(\lambda)$ is obtained from $S(\lambda)$ by induction. More precisely, $V(\lambda)$ is the quotient of the Kac representation $\overline V(\lambda)$ by its maximal sub-representation. The representation $\overline V(\lambda)$ is defined as follows. One first extend (in a trivial way) the action of ${\mathcal U}_q(\mathfrak{gl}(3)\oplus \mathfrak{gl}(1))$ to the action of an intermediate algebra and then induce this action to the whole algebra ${\mathcal U}_q(\mathfrak{gl}(3|1))$.

It follows by adjoin property that there is a non-zero map 
$$\overline V(\lambda)\to V.$$
Hence $V(\lambda)$ is a sub-quotient $V$. But they have the same character, in particular, same (total) dimension, hence are isomorphic.
\end{proof}
 \begin{lem}\label{lem_char}
The character of the representation $I(\lambda)$ constructed in the previous section is equal to the character of the highest weight irreducible representation $V(\lambda)$ of ${\mathcal U}_q(\mathfrak{gl}(3|1))$.
\end{lem}
\begin{proof}
The character of $V(\lambda)$ does not depend on $q$, hence can be computed by classical formula, for instance it is given explicitly in \cite{Dung1}. On the other hand, the character of $I(\lambda)$ can be computed directly from their construction and the compatibility of the character with exact sequences and tensor product. First, setting 
$$x_1=e^{(1,0,0|0)}, x_2=e^{(0,1,0|0)}, x_3=e^{(0,0,1|0)},y=e^{(0,0,0|1)},  $$
we have
$$\ch(I(1,0,0|0))=\chi(V)=x_1+x_2+x_3-y.$$

Using \cite[Example~I.3.22(4)]{Macdonald} we have, for $m\geq n\geq p\geq 1$,
$$\ch(I(m,n,p|0))=(x_1x_2x_3)^{p-1}(x_1+y)(x_2+y)(x_3+y)S(m-p,n-p,0)$$
where $S(m,n,p)$ is the Schur function on the variables $x_1,x_2,x_3$, associated to partition $(m,n,p)$.
Further, we have
$$\begin{array}{l}
\ch(I(m,n,0|0)
\displaystyle=\frac{(x_1+y)(x_2+y)(x_3+y)}{(x_1-x_2)(x_2-x_3)(x_1-x_3) }\times\\
\qquad\qquad\qquad\qquad\displaystyle \left(  \frac{x_2^{m+1}x_3^{n} - x_2^{n} x_3^{m+1}}{x_1 + y} 
+ \frac{x_3^{m+1}x_1^{n} - x_3^{n} x_1^{m+1}}{x_2 + y} +  \frac{x_1^{m+1}x_2^{n} - x_1^{n} x_2^{m+1}}{x_3 + y} \right), \\ \\
\ch(I(m,0,0|0)
\displaystyle=\frac{(x_1+y)(x_2+y)(x_3+y)}{(x_1-x_2)(x_2-x_3)(x_1-x_3) }\times\\
\qquad\qquad\qquad\qquad\displaystyle \left(  \frac{x_2^{m+1} - x_3^{m+1}}{x_1 + y} 
+ \frac{x_3^{m+1} -   x_1^{m+1}}{x_2 + y} +  \frac{x_1^{m+1} -   x_2^{m+1}}{x_3 + y} \right).
\end{array}$$

Since $I(1,1,1|1)$ gives the quantum super determinant, we have
$$\ch(I(1,1,1|1)=x_1x_2x_3y^{-1}.$$

 Using induction we obtain, for $k-l\neq 2$, $k\geq 2$,
 $$  \ch (\Im d_{k,l}) =\displaystyle
\frac{ (x_1+y)(x_2+y)(x_3+y)y^{k-3}}{ (x_1x_2x_3)^l} S(l ,l ,0). $$
Hence we have, according to \eqref{I5}, for $m\geq 0> p$,
\begin{equation*} 
\ch (I(m,m,p|0)=\displaystyle {(x_1+y)(x_2+y)(x_3+y)}{(x_1x_2x_3)^{p-1}}S (m-p , m-p , 0).
\end{equation*} 

Next, we have, for $i,a\geq 0$,
$$\begin{array}{l}
\ch(X_{i,a})
\displaystyle=  
\frac{(x_1+y)(x_2+y)(x_3+y)}{(x_1-x_2)(x_2-x_3)(x_1-x_3)  y} 
\left( \frac{x_1(x_2^{-a-i-1}x_3^{i+2} - x_2^{i+2}x_3^{-a-i-1})}{x_1+y}+\right.\\ \\
\displaystyle
\qquad\qquad\quad +  \left.
\frac{x_2(x_3^{-a-i-1}x_1^{i+2} - x_3^{i+2}x_1^{-a-i-1})}
{x_2+y}
 + \frac{x_3(x_1^{-a-i-1}x_2^{i+2} - x_1^{i+2}x_2^{-a-i-1})}{x_3+y}  \right)
\end{array}$$
That is, $X_{i,a}$ has the same character as the comodule $V(i+1,0,-a-i|1)$. 

Finally, we have, for $i\geq 0, k\geq 2,a+i+k\geq 0$,
$$
\ch(Y_{i,k,a})= \frac{(x_1+y)(x_2+y)(x_3+y) y^{k-3}}{ (x_1x_2x_3)^{a+i+k+1}} S(a+2i+k+2,a+i+k+1,0).
$$
That is $Y_{i,k,a}$ has the same character as $V(i+2,1,-a-i-k|3-k)$. 
This formula for the case $a+i+3\neq 0$ follows from the character formula for $d_{k,l}$ given above.

For the case $a+i+3=0$, the comodule $\Im d_{k,k-2}$ is not simple, its character can be computed by using the complex $K_2$. Indeed, have
$\Im d_{2,0}=\exterior_2.$
Using induction and the fact that the homology of $K_2$ is concentrated at the term $(3,1)$ and is $I(1,1,1|1)$ one can show that $\Im d_{k,k-2}$ has a decomposition series consisting of
$I(1,1,2-k|2-k)$ and $I(1,1,3-k|3-k)$.

 By the formulas given above one can easily check that for any dominant weight $(m,n,p|0)$
 $$I(m,n,p|0)\cong V(m,n,p|0).$$
 This finishes the proof.
\end{proof}
The following theorem is a direct consequent of the two lemmas above.

 \begin{thm} The comodules $I(\lambda)$ constructed in the previous section are simple and furnish all simple comodules  of  the Hopf super algebra $H_R$.
 \end{thm}
 
 \begin{rmk}\em  The first named author has constructed in \cite{Dung1} a full list of irreducible representations of the super group $GL(3|1)$.
There is unfortunately several misprints in that work that makes the list in fact incomplete.
The description here fulfills this gap.
 \end{rmk} 

\end{document}